\documentclass[12pt,a4paper]{amsart}
\newtheorem{theorem}{Theorem}[section]
\newtheorem{corollary}[theorem]{Corollary}
\newtheorem{lemma}[theorem]{Lemma}
\newtheorem{proposition}[theorem]{Proposition}
\newtheorem{propdef}[theorem]{Proposition/Definition}
\theoremstyle{remark}
\newtheorem{remark}[theorem]{\sc Remark}
\theoremstyle{definition}
\newtheorem{definition}[theorem]{Definition}

\theoremstyle{remark}
\newtheorem{example}[theorem]{\sc Example}

\newcommand{\gras}[1]{{\mathbb #1}}
\newcommand{\C}{\gras{C}}
\newcommand{\Z}{\gras{Z}}
\newcommand{\N}{\gras{N}}

\newcommand{\X}{\mathcal{X}}
\newcommand{\R}{\gras{R}}
\newcommand{\PP}{\gras{P}}

\begin{document}
\title[Milnor fillable contact 3-manifolds]{Milnor open books and\\ Milnor
fillable contact 3-manifolds}
\author{\sc Cl{\'e}ment Caubel}
\address{Institut de
  Math{\'e}matiques-UMR CNRS 7586, {\'e}quipe "G{\'e}om{\'e}trie et dynamique" \\case
  7012, 2, place Jussieu, 75251 Paris cedex 05, France.}
\email{caubel@math.jussieu.fr}
\author{\sc Andr{\'a}s N{\'e}methi}
\address{Renyi Institute of Mathematics, Budapest, Hungary}
\email{nemethi@renyi.hu} 
\author{\sc Patrick Popescu-Pampu}
\address{Univ. Paris 7 Denis Diderot, Inst. de
  Maths.-UMR CNRS 7586, {\'e}quipe "G{\'e}om{\'e}trie et dynamique" \\case
  7012, 2, place Jussieu, 75251 Paris cedex 05, France.}
\email{ppopescu@math.jussieu.fr}

\subjclass{32S55, 53D10, 32S25, 57R17}

\keywords{contact structures, open book decompositions, isolated
  singularities}

\date{\today}
\begin{abstract}
We say that a contact manifold $(M,\xi)$ is \textit{Milnor fillable} 
if it is contactomorphic to
the contact boundary of an isolated complex-analytic
singularity $(\X,x)$.
In this article we prove that any
3--dimensional oriented manifold 
admits at most one Milnor fillable contact structure up to contactomorphism.
The  proof is based on  \textit{Milnor open books}:
we associate with  any holomorphic function $f : (\X,x) \rightarrow
  (\C,0)$, with isolated singularity at $x$ (and any euclidian rug function  
$\rho$), an open book decomposition of $M$,  and we verify that all these 
open books  carry the contact structure $\xi$ of $(M,\xi)$ --
generalizing results of Milnor and Giroux. 
\end{abstract}

\maketitle
\pagestyle{myheadings} \markboth{{\normalsize 
C. Caubel, A. N{\'e}methi and P. Popescu-Pampu}}
{{\normalsize Milnor fillable contact 3-manifolds}}

\hfill \textit{To Bernard Teissier, for his 60-th birthday}

\section{Introduction} \label{intro}
Let $x\in \X$ be an isolated singular point of a complex analytic variety.
Choose a local embedding $(\X,x)\subset(\C^s,0)$
and intersect  $\X$ with small euclidian spheres centered at $0$.
Then the diffeomorphism type of the 
intersection is independent of the embedding and the (sufficiently small)
spheres. Hence one gets (the diffeomorphism type of) a closed oriented
manifold $M(\X)$, which is called the \emph{(abstract) boundary (or the link)
of $(\X,x)$}.

Now, for any fixed embedding $e$ and small sphere (of radius
$\sqrt{\varepsilon}$), 
the above intersection $M_{e,\varepsilon}$ 
-- being a  real 1--codimensional submanifold in the
complex manifold $\X\setminus\{x\}$ -- is naturally endowed with a complex 
hyperplane distribution $\xi_{e,\varepsilon}$, with complex multiplication
$ i|\xi_{e,\varepsilon}$. Then the  triple
$(M_{e,\varepsilon},\xi_{e,\varepsilon}, 
i|\xi_{e,\varepsilon})$ is a CR-manifold. By a result
of Scherk \cite{S 86} it is an absolute invariant: 
it determines the analytic type of the germ $(\X,x)$. On the other hand,
this CR-manifold really depends on the parameters $(e,\varepsilon)$. 

Nevertheless, one gets a well-defined intermediate invariant object by 
removing the complex multiplication $i|\xi_{e,\varepsilon}$ (but
keeping the orientation of $\xi_{e, \varepsilon}$  induced by this 
multiplication). Indeed,  Varchenko has 
showed in \cite{V 80} that the couple $(M_{e,\varepsilon},\xi_{e,
\varepsilon})$ is a contact manifold,  which, up to a contact
isotopy, only depends on the analytic type of $(\X,x)$. 
This isotopy type is called the  \emph{contact   boundary} of $(\X,x)$
and denoted by $(M(\X),\xi(\X))$.
In this way, one  associates  to any
isolated singular point of  a complex variety a contact manifold.  
(For a more general definition, see section 3.)

\begin{definition}
 Let $(M, \xi)$ be a connected closed oriented contact manifold. If
 there exists a germ $(\X,x)$ of (normal) complex analytic
 space with isolated singularity such that $(M, \xi)$ is isomorphic to
 the contact boundary $(M(\X),\xi(\X))$ of $(\X,x)$,
 then we say that $(M, \xi)$ {\em admits a Milnor filling} (or is 
 {\em Milnor fillable}). The germ $(\X,x)$ is called a {\em Milnor
 filling} of $(M, \xi)$. 
\end{definition}

We will use the same terminology for any manifold $M$ (forgetting the
contact structure) if $M$ 
is isomorphic to the abstract boundary $M(\X)$ of some
singularity $x\in\X$.

Above, the normality assumption is not restrictive: 
$(M(\X),\xi(\X))$ is isomorphic to the contact boundary
$(M(\hat \X),\xi(\hat \X))$ of the normalization  $(\hat \X,\hat x)$
of $(\X,x)$.

Some natural questions arise, the most basic one being the classification
of Milnor fillable contact manifolds. In this paper we will concentrate
on 3-dimensional oriented manifolds. In this case
the {\em existence} of a Milnor filling
is a topological property and it is completely understood:
an oriented 3--manifold $M$ is Milnor fillable if and only if it is a 
graph-manifold obtained
by plumbing along a weighted graph which has a negative definite
intersection form (see Grauert \cite{G 62}). 

\vspace{1mm}

Our main theorem establishes the {\em uniqueness} property: 

\begin{theorem} \label{mthm}
Any Milnor fillable 3-manifold admits a unique Milnor
fillable contact structure up to contactomorphism.
\end{theorem}

Here some comments are in order.

\vspace{2mm}

\noindent (i)\ 
All Milnor fillable contact structures are Stein fillable. Indeed, 
if the surface singularity $(\mathcal{S},0)$ is
smoothable, a simple application of Gray's Theorem shows that its
contact boundary  coincides with the contact boundary 
 of any Milnor fiber, which is Stein.
But even if  $(\mathcal{S},0)$ is not smoothable, its contact boundary
can still be filled with a complex manifold (e.g., with the resolution of the
singularity) and results of Bogomolov and de Oliveira (see  \cite{BO
  97}) show that the complex 
structure of the filling can be made Stein without changing the contact
boundary. In particular, Milnor fillable contact structures are tight,
by a general theorem of Gromov and Eliashberg.\\
(ii)\ 
All these manifolds contain at least one incompressible torus, except for
the lens spaces and some small Seifert spaces. Therefore, in general,
by a theorem of Colin and Honda-Kazez-Mati{\'c}
(see Colin, Giroux and Honda \cite{CGH 02}) all these
manifolds admit infinitely many different tight contact structures up
to isomorphism. In particular, theorem \ref{mthm} indicates that
Milnor fillability is a very special property of  tight contact
structures. \\ 
(iii)\ Finally notice that the above classification theorem  \ref{mthm} 
is in a big contrast with the phenomenon valid for higher dimensional
singularities.   For instance, I. Ustilovsky in \cite{U 99}  discovered
on the spheres $S^{4n+1}$, $n\geq 1$, 
infinitely many different Milnor fillable contact structures.
(Also, the authors know no criterions which  would ensure, in 
general,  Milnor fillability.)

\medskip

In section \ref{open}, we recall the work of Giroux on contact structures
and open books. The key message is that 
in dimension 3 any open book carries a unique
contact structure up to isotopy.
This description of contact structures is perfectly
adapted to the contact boundaries exploited in 
section \ref{sectionportage}. Here we start with a rather general 
definition of the contact boundary. Then we prove 
that any analytic function $f:(\X,x)\to (\C,0)$
(with an isolated singularity at $x$) defines an open book decomposition 
of the boundary $M(\X)$ --
we call them  \emph{Milnor open books}. The point is that any Milnor open book
carries the natural contact structure $\xi(\X)$   (see \ref{carry}).

We emphasize that this result is valid in {\em any} dimension and 
(we believe that) will be an essential tool 
in the further study of contact boundaries.

After all this said, the proof of theorem \ref{mthm} runs as follows
(see the end of section 4). 
For any Milnor fillable 3-manifold $M$,  using its plumbing representation,
we construct a link in it
which is isotopic to the binding of a Milnor open book for \emph{any}
Milnor filling of $M$ (see \ref{realisation}). 
Then,  using  the work of Chaves \cite{C 96} and Pichon \cite{P 01}, 
we show that all the possible Milnor open books with this 
binding  are, in fact, isomorphic  (see \ref{ubique}).

In section \ref{persp} we list some questions, 
pespectives for further studies.
\medskip

A short preliminary version of this article appeared in \cite{CP
  04}. In it, theorem \ref{mthm} \,  was proved for rational homology
spheres. In the meantime, we succeeded to replace the old proof by 
a more natural one 
and to extend the result to all Milnor fillable 3-manifolds.

\subsection*{Acknowledgements}
We are grateful to J{\'a}nos Koll{\'a}r for the simplification of the proof of
proposition \ref{realisation} and to 
Fran\c coise Michel and Anne Pichon for conversations that led to
the proof of proposition \ref{invariant}.
We would also like to thank Michel Boileau, Etienne Ghys, Emmanuel Giroux
and Bernard Teissier for having kindly answered our questions.

\section{Contact structures and open books}
\label{open}

Let $M$ be an oriented $(2n-1)$-dimensional manifold. A 
\emph{contact structure}
on $M$ is a hyperplane distribution $\xi$ in $TM$ given by a global
1-form $\alpha$
such that $\alpha\wedge(d\alpha)^{\wedge (n-1)}\neq0$. 
We say that the pair $(M,\xi)$ is a  \emph{contact manifold} and 
$\alpha$ a \emph{contact form}. 
The form $\alpha$ is called \textit{positive} if
$\alpha\wedge(d\alpha)^{\wedge (n-1)}$ defines the chosen orientation
of $M$.  If $n$ is even, then the orientation defined
by $\alpha\wedge(d\alpha)^{\wedge (n-1)}$ does not depend on the
choice of the defining form $\alpha$, hence 
one can speak about  \textit{positive contact structures}.

In the sequel we consider only \textit{oriented} contact structures,
i.e. $\xi$ denotes a field of oriented hyperplanes. 

If $\alpha$ is a contact form on $M$, the condition $\alpha \wedge
(d\alpha)^{\wedge (n-1)} \neq 0$ implies that $d\alpha|_{\xi}$ is a
symplectic form. This shows that $\ker(d \alpha)$ is a 1-dimensional
vector subspace of $TM$, transversal to $\xi$. Therefore, there exists a
unique vector field $R$ on $M$ such that $ d\alpha(R, \cdot)=0$ and 
$ \alpha(R)=1$. 
It is called the \textit{Reeb vector field} associated to $\alpha$.

Two contact structures $\xi$ and $\xi'$ on $M$ are \emph{isotopic}
(resp. \emph{isomorphic} or \emph{contactomorphic})
if there is an isotopy (resp. a diffeomorphism) of $M$ which sends
$\xi$ on $\xi'$.

For more about contact structures, see e.g. 
Eliashberg and Misha\-chev's book \cite{EM 02}. 
\medskip

In the study of  contact manifolds one of the main tools is provided by 
 \textit{open books carrying  contact structures}. 

\begin{definition}
  An {\em open book} with {\em binding} $N$ in a manifold $M$ is
  a couple $(N, \theta)$,
  where $N$ is a (not necessarily connected) 2--codimensional closed
  oriented submanifold of $M$ with trivial normal bundle
and $\theta:M\setminus N \rightarrow \mathbf{S}^{1}$ is a
  smooth
  fibration which in a neighborhood $N \times \mathbf{D}^{2}$ of $N$
  coincides with the angular coordinate.
The fibers of $\theta$ are called the {\em pages} of the open book.

We say that the 
open books $(N, \theta)$ and $(N', \theta')$ in the manifolds $M$,
respectively $M'$, are {\em isomorphic} if there exists
a diffeomorphism $\phi: (M,N) \rightarrow (M', N')$ which preserves
the orientations and carries the fibers
of $\theta$ into the fibers of $\theta'$.
\end{definition}

Notice that $d\theta$ induces natural co-orientations
on the binding and the pages of the open book. Thus, any fixed orientation of
$M$ induces a natural orientation on $N$. 
If $N$ itself is oriented \emph{a priori}, then we say that the open
book is \textit{compatible with the orientations of $M$ and $N$} if the
two orientations of $N$ coincide.

The next criterion allows us to recognize open
books. Its proof is straightforward.

\begin{lemma} \label{crit1}
   Let $M$ be an oriented closed manifold and let $\phi : M
   \rightarrow \C$ be a differentiable function. If there
   exists a number $\eta > 0$ such that

$ \bullet \: d(\arg \: \phi) \neq 0$ if $|\phi | \geq \eta$, and 

$ \bullet \: d(\phi) \neq 0$ if $| \phi | \leq \eta$\\
then $(\phi^{-1}(0), \arg \: \phi)$ is an open book in $M$.
\hfill\qed
\end{lemma}

\begin{remark}  \label{argument}
  With the exception of \ref{lambdarg}, where one chooses
  the branch $\arg\in (-\pi, \pi]$, the argument of a non-zero
  complex number is regarded  as an element of $\R/2 \pi \Z
  \simeq \mathbf{S}^1$.
\end{remark}

Next we recall the relationship between contact structures and open books.

\begin{definition} (Giroux \cite{G 02})
 We say that a positive contact structure $\xi$ on an oriented manifold
$M$ is {\em carried by an open book}
 $(N, \theta)$ if it admits a defining contact form $\alpha$ which
 verifies the following:

  $\bullet$ $\alpha$ induces a positive contact structure on $N$;

  $\bullet$ $d\alpha$ induces a positive symplectic structure on each fiber of
        $\theta$.\\
If a contact form $\alpha$ satisfies these conditions we say that 
it is {\em adapted} to $(N,\theta)$.
\end{definition}

One has the following  criterion for an open book to carry a contact
structure: 

\begin{lemma} \label{crit2} (Giroux \cite{ G 03}) Let $\alpha$ be a
  contact form on the manifold $M$. Suppose that there exists an open
  book $(N, \theta)$ in $M$ and a neighborhood $V=N \times
  \mathbf{D}^{2}$ of $N$ such that:

$\bullet$ $\theta$ is the normal angular coordinate in $V$;

$\bullet$ $\alpha$ induces a positive contact structure on each
submanifold $N \times \{*\}$ in $V$;

$\bullet$ $d\alpha$ induces a positive symplectic form on each fiber
of $\theta$ in $M \setminus \: \mathrm{Int} V$.\\
Then the open book $(N, \theta)$ carries the contact structure $\xi=
\ker \: \alpha$.\hfill\qed
\end{lemma}

The relevance of this notion for the study of 3-dimensional contact
manifolds lies in the following result:

\begin{theorem}
\label{giroux} (Giroux \cite{G 02}, \cite{G 03})
On a closed oriented 3-manifold, two positive contact
 structures carried by  the same
 open book are isotopic.\hfill\qed
\end{theorem}

In particular, in order
 to show that two contact structures on a given 3-manifold are
 isomorphic, it is enough to show that they are carried by isomorphic
 open books.

\begin{remark} Giroux also proved  that the same statement holds in all
 dimensions if one further asks
that the two contact structures induce the same symplectic structure
on the pages
up to isotopy and completion (see again \cite{G 02} and \cite{G
  03}). He also proved in collaboration with Mohsen that any contact
structure is carried by an open book. Hence, in fact in any dimension, 
one can translate statements of contact geometry into properties of 
open books. 
\end{remark}

\section{Contact boundaries and Milnor open books}
\label{sectionportage}

Let $(\X,x)$ be an irreducible  germ of a complex analytic
 space with isolated singularity. 
Sometimes, we will denote by  $\X$  a sufficiently small
representative of this germ.
 Let $m_{\X,x}\subset \mathcal{O}_{\X,x}$ be the ideal of
 germs of holomorphic 
  functions on $(\X,x)$ vanishing at $x$.

\subsection*{A. The contact boundary associated with a holomorphic immersion.}
Write $\X^*$ for the complex manifold $\X\setminus\{x\}$. 
Let  $J: T\X^* \rightarrow T\X^*$ be 
the operator of fiberwise multiplication by $i$, when
$T\X^*$ is seen as a real vector bundle. We will also denote
it by $i\cdot$, when no confusion is possible.
Set  $d^{\mathbf c}:= J^* \circ d$,
i.e. $d^{\mathbf c}F = dF \circ J$ for any differentiable
function $F:\X^{*}\to\R$. Then $d^{\mathbf c}=i (\partial
-\bar{\partial})$. 

A real function $F$ on $\X^*$ is called \textit{strictly
  pluri-subharmonic} (\textit{spsh}) if and only if $-dd^{\mathbf c}(F)>0$,
  that is if $-dd^{\mathbf c}(F)(v, Jv)>0$ for any non-zero tangent
  vector $v$ of  $\X^*$.

For any  $\phi_1,..., \phi_{N} \in m_{\X, x}$ consider the 
holomorphic map $\Phi: (\X,x) \rightarrow (\C^N, 0)$ with components $\phi_i$,
and the  real analytic function
  $$ \rho:= \sum_{k=1}^N |\phi_k|^2: (\X,x)\to (\R,0).$$
For each $\varepsilon >0$, define 
  $$M_{\rho, \varepsilon}:= \rho^{-1}(\varepsilon).$$
Clearly, $M_{\rho, \varepsilon}$ is a 
smooth compact manifold for $\varepsilon >0$ sufficiently small if and only if
$\Phi$ is a \textit{finite} analytic morphism. In the sequel we will
assume that this fact holds. 

On $\X^*$ we consider the following natural objects associated
with $\rho$:
   \begin{eqnarray*}
          \alpha&:=& -d^{\mathbf c}\rho\\
          \omega&:=& d\alpha=-dd^{\mathbf c}\rho\\
          g(u,v)&:=& \omega(u, Jv),\quad \forall \: u,v \in T\X^*\\
          h&:=& g + i \omega.
      \end{eqnarray*}

\noindent Then, on $\X^*$ define:
    $$\xi_{\rho}:= \ker(d\rho) \cap \ker (d^{\mathbf c}\rho).$$
It is a field of complex tangent hyperplanes of the real tangent
bundle of $\X^*$ with its canonical almost complex structure.  
Moreover, 
it is tangent to the levels $M_{\rho, \varepsilon}$ of $\rho$. In fact:
 $$  \xi_{\rho,\varepsilon}:=\xi_{\rho}|_{M_{\rho,\varepsilon}}= \ker
     (\alpha|_{M_{\rho, \varepsilon}}). $$

\begin{lemma} \label{condcont}
The following conditions are equivalent:

 1) The pair $(M_{\rho, \varepsilon}, \xi_{\rho, \varepsilon})$
  is a contact manifold for $\varepsilon$ sufficiently small.
 
 2) The morphism $\Phi$ is an immersion of $\X^*$ into $\C^N$. 

 3) The function $\rho$ is spsh.
\end{lemma}

  \begin{proof}
      An easy computation shows that
    $$-dd^{\mathbf c}\rho(v,Jw)= 2 \sum_{k=1}^{N}\det
         \left( \begin{array}{cc}
                    d\phi_k(v) & -\: d\phi_k(w)\\
\\[-3mm]
                   \overline{d\phi_k(v)} & \: \overline{d\phi_k(w)}
                 \end{array}  \right)$$
for any tangent vector fields $v,w$ of $\X^*$. 
  This shows that $\ker \, d \Phi = \ker (-dd^{\mathbf c}\rho)$ 
and  $-dd^{\mathbf c}\rho \geq 0$. The lemma is an immediate consequence of
  this.
  \end{proof}

 From now on we fix $\Phi$, and we assume that it induces  an 
\textit{immersion of} 
$\X^*$ into $\C^N$; we will briefly say that it is \textit{a holomorphic
  immersion}.  In such a  case, we say that 
the function $\rho:\X^*\to\R$ defined above is an 
\emph{euclidian rug function}.
In the view of the previous lemma, the objects associated with 
an euclidian rug function $\rho$
have the following properties: $\alpha$ -- restricted to the levels
$M_{\rho,\varepsilon}$ -- is a contact form, $g$ is a
riemannian metric, $\omega$ is a symplectic form compatible with
the complex structure on $\X^*$ and $h$ is a hermitian metric.   
(Notice that not all the real analytic rug functions, e.g. as in 
\cite{L 84}, have all these properties, although they are perfectly good 
to identify the boundary manifold $M_{\rho,\varepsilon}$.)

\medskip

Now, using Gray's theorem (see \cite{EM 02} for instance), it is easy to prove:

\begin{propdef}\label{defcont}
The pair $(M_{\rho,\varepsilon},\xi_{\rho,\varepsilon})$ is a positive contact
manifold, whose contact isotopy type  does not 
depend on the choice of the holomorphic immersion $\Phi$ and of
$\varepsilon>0$ sufficiently small. This isotopy type  is called  
the {\em contact boundary
of $(\X,x)$} and  denoted by $(M(\X),\xi(\X))$.\hfill\qed
\end{propdef}

This generalizes Varchenko's result  \cite{V 80} which  corresponds 
to the case when  $\Phi$ is an embedding. 
Considering arbitrary holomorphic immersions and their associated rug
functions (instead of embeddings and euclidean spheres) has many
advantages: we not only increase the possibilities to realize the isotopy
type of the natural contact structure, but it also gives possibility to 
compare the corresponding contact structures under finite coverings.

\begin{proposition} \label{quot}
  Suppose that $(\X,x)$ is normal. Let $G$ be a finite group of
  analytic automorphisms of $(\X,x)$, acting freely and effectively
  on
  $\X^*$. Denote by $(\mathcal{Y},y)$ the quotient of $(\X,x)$ by 
   $G$. Then the isotopy type of 
  $(M(\X),\xi(\X))$ admits  a representative 
  on which $G$ acts via contactomorphisms, and the quotient  
  is a representative of   $(M(\mathcal{Y}), \xi(\mathcal{Y}))$.
\end{proposition}

\begin{proof}
Let $\gamma:(\X,x) \rightarrow (\mathcal{Y},y)$ be the quotient
map. Choose an embedding $\Psi: (\mathcal{Y},y) \rightarrow \C^M$ and
define $\nu:= \sum |\psi_k|^2$, $\rho:= \nu \circ \gamma$, $\Phi:=
\Psi\circ \gamma$. Then  compare the contact structures associated with these 
rug functions. 
\end{proof}

\subsection*{B. The Milnor open book associated with a holomorphic function.}
Fix an  $f \in m_{\X,x}$ which defines an isolated singularity at
$x$. Set 
$$N_{\rho,\varepsilon}(f):=M_{\rho,\varepsilon}\cap f^{-1}(0).$$
For $\varepsilon >0$ sufficiently small $N_{\rho,\varepsilon}(f)$ 
is smooth and naturally oriented.
The argument of $f$ restricted to $M_{\rho,\varepsilon}\setminus
N_{\rho,\varepsilon}(f)$ 
gives a well-defined function
$$\theta_\varepsilon(f):=\mbox{arg }f:M_{\rho,\varepsilon}\setminus
N_{\rho,\varepsilon}(f)\to \mathbf{S}^1.$$
We then have the following generalization of Milnor's Fibration
Theorem (see \cite{M 68}):

\begin{propdef}\label{defmilnor}
For $\varepsilon>0$ sufficiently small, the pair
$(N_{\rho,\varepsilon}(f),\theta_\varepsilon(f))$
is an open book in the boundary $M_{\rho,\varepsilon}$ which is
compatible with the orientations.
Furthermore, its isotopy type  does not
depend on the choice of $\varepsilon>0$ nor on
the choice of holomorphic immersion $\Phi$.
It  is called the {\em Milnor open book  of $f$} and 
denoted by $(N(f),\theta(f))$. The pair $(M(\X),N(f))$ is called the
{\em link} of $f$.
\end{propdef}

\noindent 
This statement splits into two parts: a fibration and an
invariance result. The first one, for $\Phi$ an embedding, 
appears  in \cite{H 71}, Satz 1.6.
That proof extends to our new situation once a key fact is verified.
Since this fact is used by us in other places as well, we provide its 
complete proof in proposition \ref{lambdarg}, after having stated
two preliminary lemmas. 
For all the other details, we refer to \cite{H 71}.
The invariance statement can be proved similarly using the classical tools
of local analytic singularities: we leave the verification 
to the interested reader.

\begin{remark}
(a) In fact,  a  fibration result has been proved
in a more general context by Durfee in
\cite{D 83}, but without specifying the fibration map: he actually proves
that the complement of  $f^{-1}(0)$ in the boundary of any analytic
neighborhood of a (non necessarily isolated) singularity $x$ in
$\X$ is  the total space of a fibration
(with no more precision about the projection map) 
whose fiber is homeomorphic to the intersection
of a smooth fiber of $f$ with this neighborhood.
The above proposition \ref{defmilnor}  shows that, 
in the neighborhoods defined by euclidian rug functions,
this fibration can be defined by the argument of $f$.
A corresponding statement in the general case is not guaranteed: 
we use indeed the peculiar form of $\rho$ in the computations in 
\ref{lambdarg}.
We also emphasize that if one wishes to verify the compatibility of an
open book with a contact structure, then one needs very precise information
about $\theta$ and about (a ``well chosen'') contact form $\alpha$. 
This actually explains why we need  the fibration  to be 
given by the argument of $f$.\\
(b) Fix the boundary $M$ of an analytic germ $(\X,x)$ as above.
In general it is  extremely difficult to verify that an open book
$(N,\theta)$ in $M$ (or a link $N$ in $M$) is isotopic with a 
Milnor open book (or link) of a function  germ $f$ on $(\X,x)$. 
Even in the surface case it can happen that some 
open book is determined by a function germ for some analytic structure 
of $(\X,x)$, but the same fact is not true if one modifies the 
analytic structure of $(\X,x)$  (see \cite{N 99}, (2.15)).
\end{remark}

For any real function $F$ defined on  $\X^*$, its gradient $\nabla F$
will be taken with respect to the riemannian metric $g$. If $\phi \in
m_{\X,x}$, we also denote by $\nabla \phi$ its gradient with
respect to the hermitian metric $h$, that is $d\phi = h(\nabla \phi,
\cdot)$ ($\phi$ being holomorphic, this field is well-defined). 

\vspace{1mm}

The proof of the following lemma is straightforward:
\begin{lemma}\label{diff}
If $\phi \in m_{\X,x}$, one has:\\
1) $\nabla|\phi|^2=2\phi\nabla \phi$;\\
2) In $\X\setminus \phi^{-1}(0)$, $d \arg \phi = \mathrm{Im}\left(\dfrac{d
  \phi}{\phi}\right)$ and 
$\nabla \arg \phi=i\dfrac{\nabla\phi}{\bar\phi}$.\hfill\qed
\end{lemma}

\begin{lemma} \label{arc}
  Let $p: (\R_+, 0) \rightarrow ( \X, x)$ be a non-constant 
real analytic arc, and set $\dot{p}:= \frac{dp}{dt}$. Then:
      $$\lim_{t\rightarrow 0} \: h\left(\frac{\nabla \rho}{\|\nabla
            \rho\|}(p(t)),   
            \frac{\dot{p}(t)}{\| \dot{p} (t)\|}\right) =1.$$
\end{lemma}

\begin{proof}
   Let $h_0$ be the canonical hermitian form on $\C^N$. If
   $\rho_0= \sum_{k=1}^{N}|z_k|^2$, then its imaginary part is
   $\omega_0= -dd^{\mathbf c}\rho_0$. As $\rho=
   (\rho_0|_{\Phi(\X)})\circ\Phi $, we get: 
     $$\begin{array}{l}
          h= \Phi^* (h_0|_{\Phi(\X)})\\
          \Phi_*(\nabla \rho)= \nabla_0 (\rho_0|_{\Phi(\X)})\\
          \Phi_*(\dot{p})= \dot{q},
        \end{array}$$
   where $q:= \Phi\circ p$ and $\nabla_0$ is the gradient
   with respect to the riemannian 
metric $g_0:=\mathrm{Re}\,h_0|_{\Phi(\X)}$.
Hence, we get the following equalities of functions of $t$ (where we write
briefly $\nabla \rho$ instead of $(\nabla \rho)(p(t))$):
 $$
       h\left(\frac{\nabla \rho}{\|\nabla \rho\|}, 
            \frac{\dot{p}}{\| \dot{p} \|}\right)=
     \Phi^* h_0\left(\frac{\nabla \rho}{\|\nabla \rho\|}, 
            \frac{\dot{p}}{\| \dot{p} \|}\right)= 
      h_0\left(\Phi_*\left(\frac{\nabla \rho}{\|\nabla \rho\|}\right), 
            \Phi_*\left(\frac{\dot{p}}{\| \dot{p} \|}\right)\right) $$
  $$
       = h_0\left(\frac{\nabla_0 (\rho_0|_{\Phi(\X)})}{\|\nabla_0(
            \rho_0|_{\Phi(\X)})\|},  
            \frac{\dot{q}}{\| \dot{q} \|}\right). $$
Take now 
a semi-analytic neighborhood $U(p)$ of the image of $p$ in $\X^*$ which
is embedded in $\C^N$ by $\Phi$. Then the pair $(\Phi(U(p)), \{0\})$
verifies Whitney's condition (a) at $0$ (see e.g. \cite{L 84}),
which shows that the angle
between the vector $q(t) \in \C^N$ and the tangent space to
$\Phi(U(p))$ at the point $q(t)$ converges to $0$ when $t \rightarrow
0$. This implies that: 
   $$\frac{\nabla_0 (\rho_0|_{\Phi(\X)})}{\|\nabla_0(
            \rho_0|_{\Phi(\X)})\|} (q(t))= \frac{\nabla_0\rho_0}{\|\nabla_0
            \rho_0\|}(q(t)) +o(1)=  \frac{q}{\| q \|}(t) + o(1),$$ 
where $\nabla_0\rho_0$ denotes the gradient on the ambient space
            $\C^N$. 
From the previous computation we get:
 \begin{equation} \label{asympt}
    h\left(\frac{\nabla \rho}{\|\nabla \rho\|}, 
            \frac{\dot{p}}{\| \dot{p} \|}\right)= h_0\left( \frac{q}{\| q \|},
            \frac{\dot{q}}{\| \dot{q} \|}\right) + o(1). 
 \end{equation}
Let $q(t)= q_0 t^Q + o(t^Q)$ be the beginning of the series
expansion of $q(t)$, where $q_0\in \C^N-\{0\}$ and $Q \in \N^*$. Then 
 $\dot{q}(t)= Qq_0t^{Q-1} + o(t^Q)$, which implies:
  $$ \frac{q}{\| q \|}= \frac{q_0}{\| q_0 \|} +o(1), \quad
       \frac{\dot{q}}{\| \dot{q} \|} = \frac{q_0}{\| q_0 \|} +o(1).$$
Using equation (\ref{asympt}), the lemma is proved.
\end{proof}

The next proposition -- which is the key step in the proof
of  \ref{defmilnor} --  generalizes lemma 4.3 of
\cite{M 68} to the case of singular ambient spaces 
 and {\em immersions}  $\Phi$ (which are not necessarily embeddings). 
Our proof runs rather similarly, with 
the difference that our computations are intrinsic, they do not depend
on any choice of local coordinates. 

\begin{proposition} \label{lambdarg}
 Let $\phi \in m_{\X, x}$.   
 For any $\theta_0 \in ( 0, \frac{\pi}{2})$, there exists a
  neighborhood $U_{\theta_0}$ of $0$ in $\X$ such that inside
  $U_{\theta_0}\setminus \phi^{-1}(0)$ the following 
  implication  holds:
     \begin{equation*} 
     \nabla \arg\phi =i\lambda \nabla \rho \mbox{ with } \lambda \in
     \C^* \: \Longrightarrow  | \arg \lambda | < \theta_0.
     \tag{$*$} \end{equation*} 
\end{proposition}

\begin{proof}
  Suppose, by contradiction,  that ($*$)  does not hold. 
  Then, by the Curve Selection Lemma (see \cite{L 84}),
  there exists an 
  analytic arc $p: (\R_+,0) \rightarrow (\X,x)$ such that along it we
  have the equality:
     \begin{equation} \label{proport}
         \nabla (\arg\phi) (p(t))=i\lambda(t) \nabla \rho (p(t))
     \end{equation}
with $|\arg  \lambda(t)| \geq \theta_0$ for any sufficiently small $t$.

As $\nabla \arg\phi = i\nabla \phi/\bar\phi$  by lemma \ref{diff}, part 2, 
we get:
  $$ \frac{d}{dt} \phi(p(t)) =h\left(\nabla \phi, \frac{dp}{dt}\right)=
         h(-i\bar\phi\nabla \arg\phi, \dot{p})= 
           h(\bar\phi\lambda \nabla \rho, \dot{p})=
            \phi\bar\lambda
              h(\nabla \rho, \dot{p}).$$
This implies:
  \begin{equation} \label{lambda}
      \lambda(t)= \left(\overline{\phi(p(t))}\right)^{-1}\cdot
            \frac{d}{dt} \overline{\phi(p(t))} 
            \cdot \|\nabla \rho\|^{-1}\cdot \|\dot{p}\|^{-1}\cdot 
                h\left(\frac{\nabla \rho}{\|\nabla \rho\|}, 
            \frac{\dot{p}}{\| \dot{p} \|}\right)^{-1}.
  \end{equation}
The equation (\ref{lambda}) shows that the function $\lambda(t)$ has a
Laurent expansion of  type:
   \begin{equation} \label{laurlambda}
         \lambda(t)=l t^L + o(t^L), \: l\neq 0, \: L \in \Z.
   \end{equation}
Consider now the other two expansions  as well:
  $$ \begin{array}{c}
        \phi(p(t))= at^A +o(t^A), \: a\neq 0, \: A \in \N^*\\[1mm]
        \|\nabla \rho\|^{-1}\cdot \|\dot{p}\|^{-1} = bt^B + o(t^B), \:
        b>0, \: B \in \Z.
      \end{array}$$
From the first one we deduce:
  $$\frac{d}{dt} \overline{\phi(p(t))}= A\overline{a}t^{A-1} +
  o(t^{A-1}).$$
Combining them with lemma \ref{arc} and with equation (\ref{lambda}),
  we get:
   \begin{equation} \label{exprl}
     l= 
    A\cdot b >0
   \end{equation}
Since (\ref{laurlambda}) and (\ref{exprl}) contradict the
hypothesis $|\arg \: \lambda(t)| \geq \theta_0$, the proposition is proved.
\end{proof}

\subsection*{C. The Milnor open books carry the natural contact
  structure on the boundary.} 
 Let us summarize: we have associated to any isolated singularity 
$(\X,x)$ a well-defined contact structure on
its boundary $M_{\rho,\varepsilon}$, and to any function $f\in
 m_{\X,x}$ with an isolated singularity
an open book $(N_{\rho,\varepsilon}(f),\theta_\varepsilon(f))$ on
 $M_{\rho,\varepsilon}$.  
These two objects are naturally related:

\begin{theorem} \label{carry}
Let $(\X,x)$ be a complex analytic variety having
an isolated singularity at $x$, and let $f:(\X,x)\to(\C,0)$
be a holomorphic function having also an isolated singularity.
Then the Milnor open book $(N_{\rho,\varepsilon}(f),\theta_\varepsilon(f))$
of $f$ carries the natural contact structure
on the boundary $M_{\rho,\varepsilon}$  of $\X$.
\end{theorem}

\begin{remark}
(a) Theorem \ref{carry} stengthens a result of Giroux (see
\cite{G 03} for more details) and generalizes it to a singular ambient space:
Giroux's  original proof is  valid only up to isotopy,  the contact
boundary $M_{\rho,\varepsilon}$ being replaced there by one of its
deformations  
-- a level of the function $\rho_c:=\rho+c|f|^2$, for $c\gg 1$.\\
(b) The particular case of the function $z_0\in m_{\X_k,0}$ with
$\X_k=\{z_0^k+z_1^2+\ldots +z_n^2=0\}\subset\C^{n+1}$ has been studied
by Van Koert and Niederkr{\"u}ger in
\cite{vKN 04}, in relation with Ustilovsky's spheres.\\
(c) Theorem \ref{carry} has the following consequence. Let us fix the 
analytic germ $(\X,x)$. Then all the open books 
associated with all the possible  holomorphic function germs $f$
(with isolated singularity at $x$)  
determine (up to isotopy) the same contact structure. 
Notice also 
that function germs $f$ with isolated singularity always exist.
Indeed, once an embedding of $(\X,x)$ into $(\C^N,0)$
    is chosen, it is enough to take the restriction to $\X$ of a
    linear form whose kernel is not a limit of tangent hyperplanes to
   $\X\setminus\{x\}$ (see L{\^e}, Teissier \cite{LT 88} for details).
\end{remark}

We start the proof of \ref{carry}  with some lemmas.

Fix an euclidian rug function $\rho:\X\to\R$. For $\varepsilon$ sufficiently
small, the 1-form
$\alpha=-d^{\mathbf c}\rho$ defines the natural contact structure
on the smooth level $M_{\rho,\varepsilon}=\rho^{-1}(\varepsilon)$.
Denote $R\in\Gamma(M_{\rho,\varepsilon},TM_{\rho,\varepsilon})$
its Reeb vector field.

\begin{lemma}\label{reeb}
The Reeb vector field $R$ of $\alpha$ is given by
$R= i\nabla\rho/\|\nabla\rho\|^2$. 
Moreover, the contact distribution $\xi_{\rho,\varepsilon}$ on
$M_{\rho,\varepsilon}$ is exactly the orthogonal 
complement of $\C.R=\C.\nabla\rho$ in $T\X^*|_{M_{\rho,\varepsilon}}$
with respect to the hermitian 
form $h$ associated with $\rho$.
\end{lemma}

\begin{proof}
Since $R$ is a generator of the kernel of $\omega$ restricted to
$TM_{\rho,\varepsilon}$, on $T\X^*|_{M_{\rho,\varepsilon}}$  we have 
 $\iota_R\omega=kd\rho$ for some 
$k\in\mathcal{C}^\infty(M_{\rho,\varepsilon},\R)$.
This shows that on $T\X^*|_{M_{\rho,\varepsilon}}$ one has 
$$\iota_R\omega=g(k\nabla\rho,\cdot)=\omega(-ki\nabla\rho,\cdot).$$
Since $\omega$ is non-degenerate on
$T\X^*|_{M_{\rho,\varepsilon}}$, this shows that 
$R=-ki\nabla\rho$. Hence  
$$1=\alpha(R)=-d\rho(iR)=-d\rho(k\nabla\rho)=-k\|\nabla\rho\|^2,$$
which proves the first statement.
For the second statement, it suffices to notice that
$$h(R,v)=\omega(R,i.v)+i\omega(R,v)=0$$
for any section $v$ of $\xi_{\rho,\varepsilon}$ (here we use that
$i.v$ is also a section of $\xi_{\rho,\varepsilon}$). 
\end{proof}

\begin{lemma}\label{formula}
Fix $c>0$, and put $\alpha_c:=e^{-c|f|^2}\alpha$. If $R_c$ denotes the Reeb
vector field of $\alpha_c$,  on $M_{\rho,\varepsilon}\setminus
N_{\rho,\varepsilon}(f)$ one has: 
$$d\theta(R_c)=e^{c|f|^2}\left(d\theta(R)+2c|f|^2\|\mathrm{pr}_\xi
\nabla\theta\|^2\right),$$ 
where $\mathrm{pr}_\xi:T\X^*|_{M_{\rho,\varepsilon}}\to\xi$ denotes
the projection parallel to $\C.R$. 
\end{lemma}
\begin{proof}
Put $H:=e ^{-c|f|^2}$. 
Put also $R_c:=k(R+S_c)$, where
$S_c\in\Gamma(M_{\rho,\varepsilon},\xi_{\rho,\varepsilon})$ and  
$k\in\mathcal{C}^\infty(M_{\rho,\varepsilon},\R)$. In fact,
$$1=\alpha_c(R_c)=H\alpha(k(R+S_c))=kH\alpha(R)=kH,$$
hence  one has $k=1/H$. Now,
$$d\alpha_c=dH\wedge\alpha+Hd\alpha$$
which, when applied to $R_c=(1/H)(R+S_c)$ and restricted to
$\xi:=\xi_{\rho,\varepsilon}=\ker\alpha$ gives 
$$\left(\iota_{S_c}d\alpha\right)|_\xi=\left.\frac{dH}{H}\right|_\xi=
-c\,d|f|^2|_\xi.$$   
But on $T\X^*$, lemma \ref{diff} implies that
$$d|f|^2=g(\nabla|f|^2,\cdot)=\omega(-i\nabla|f|^2,\cdot)=
\omega(-2|f|^2\nabla\theta,\cdot).$$
In particular,
$\left(\iota_{S_c}d\alpha\right)|_\xi=\left.\omega(2c|f|^2\nabla\theta,
  \cdot)\right|_\xi$.  
But $d\alpha|_\xi=\omega|_\xi$ is non-degenerate and
$\iota_v\omega|_\xi=\iota_{\mathrm{pr}_\xi(v)}\omega|_\xi$ 
 for any $v\in
T\X^*|_{M_{\rho,\varepsilon}}$. Hence, we get 
$$S_c=\mathrm{pr}_\xi(2c|f|^2\nabla\theta),$$
which shows that 
$$d\theta(S_c)=2c|f|^2d\theta(\mathrm{pr}_\xi\nabla\theta)=
2c|f|^2g(\nabla\theta,\mathrm{pr}_\xi\nabla\theta)=2c|f|^2\|\mathrm{pr}_\xi
\nabla\theta\|^2,$$  
the last equality being a consequence of the second statement of the
preceding lemma \ref{reeb}. 
Since $d\theta(R_c)=d\theta\left(e^{c|f|^2}(R+S_c)\right)$, we are done.
\end{proof}

\begin{proof}[Proof of Theorem \ref{carry}]
Fix a sufficiently small representative of $\X$ so that $\mathrm{Re}\,\lambda
>0$ on $\X\setminus f^{-1}(0)$ whenever  a relation of the form 
$\nabla \theta =i \lambda \nabla \rho$ holds (see
\ref{lambdarg}).
Consider also $\varepsilon >0$ sufficiently small. 
Now let $\eta>0$ be sufficiently  small to ensure that all the fibers
$f^{-1}(t)\subset\X$ 
cut $M_{\rho,\varepsilon}$ transversally for $|t|^2\leq\eta$. Denote
$$V_\eta:=
M_{\rho,\varepsilon}\cap\{|f|^2\leq\eta\}$$
 the corresponding tubular neighborhood of the binding
$N_{\rho,\varepsilon}(f)$ in $M_{\rho,\varepsilon}$. Then clearly
$\theta$ is a normal angular coordinate 
on this neighborhood $V_\eta\simeq N_{\rho,\varepsilon}(f)\times
\mathbf{D}^2_{\sqrt{\eta}}$. 
Moreover, each submanifold $(f^{-1}(t)\cap
M_{\rho,\varepsilon})\subset V_\eta$, 
being a level of the strictly plurisubharmonic function $\rho$ on the
complex manifold 
$f^{-1}(t)\setminus\{0\}$, is a positive contact submanifold of
$M_{\rho,\varepsilon}$. 
Thus lemma \ref{crit2} will imply theorem \ref{carry} if we can find a
convenient $c>0$ 
such that $d\alpha_c$ induces a positive symplectic form on each fiber
of $\theta$ in 
$M_{\rho,\varepsilon}\setminus\mathrm{Int}V_\eta$. But this is exactly
equivalent to the inequality  
$d\theta(R_c)>0$. 

Now, there is a $m>0$ such that $d\theta(R)\geq -m$ on the compact set
$M_{\rho,\varepsilon}\setminus\mathrm{Int}V_\eta$.
Put 
$$Z_\varepsilon:=(M_{\rho,\varepsilon}\setminus\mathrm{Int}V_\eta)\cap
\{d\theta(R)\leq 0\} \ \ \mbox{and } \ \ 
k:=\min_{Z_\varepsilon}(|f|^2\|\mathrm{pr}_\xi\nabla\theta\|^2).$$
Since $Z_\varepsilon$ is compact, $k$ is well defined.
If $k>0$, then for $c:=m/k$ we will always have $d\theta(R_c)>0$ on
$M_{\rho,\varepsilon}\setminus\mathrm{Int}V_\eta$
by lemma \ref{formula}, and the theorem is proved.

Next, assume $k=0$. This means that there is a $p\in Z_\varepsilon$ so that
$\mathrm{pr}_{\xi_p}\nabla\theta(p)=0$. 
But this implies that 
$\nabla\theta(p)=i\lambda\nabla\rho(p)$
for some  $\lambda \in \C$. Our initial choice of the representative $\X$
(see \ref{lambdarg}) implies that 
$\mathrm{Re}\, \lambda >0$. But then 
\begin{eqnarray*}
d\theta(p)(R(p))&=&
g(\nabla\theta(p),i\nabla\rho(p)/\|\nabla\rho(p)\|^2)\\
&=&\frac{1}{\|\nabla\rho(p)\|^2}\,g(i\lambda\nabla\rho(p),i\nabla\rho(p))\\
&=&\mathrm{Re}\,\lambda>0,
\end{eqnarray*} which is impossible since $p\in Z_\varepsilon$.
\end{proof}

\section{Construction of an ubiquitous Milnor open book}
\label{construction}

In this section we will restrict ourselves to normal surface singularities
(see \cite{N 99} for further references) and we construct 
for any given 3-manifold $M$  -- 
which is the boundary of some Milnor filling -- 
an open book decomposition,  which is isomorphic
to a Milnor open book for \emph{any} Milnor filling of $M$.

\subsection*{A. A sufficient condition to be the exceptional 
part of the divisor of a function.}
We start with some notations. 
Let $(\mathcal{S},0)$ be the  germ of a normal complex analytic surface 
singularity. Fix a {\em good} resolution 
$p:(\Sigma, E) \rightarrow (\mathcal{S},0)$. Namely, $\Sigma$ is
smooth, $p$ is proper and realizes an isomorphism over
$\mathcal{S}-\{0\}$ and finally the set-theoretical fibre  
$E:=p^{-1}(0)$ is a normal crossing divisor in $\Sigma$ having smooth
irreducible components $E_{1},\ldots,E_{r}$. For each $i$, we denote by
$g_i$ the genus, respectively by  $v_i:=E_i\cdot(E-E_i)$ the 
valency of $E_i$ seen as a vertex of the dual graph of $E$. 
In general, we prefer to  fix a Stein representative ${\mathcal S}$ of
$(\mathcal{S},0)$ and to set $\Sigma=p^{-1}(\mathcal{S})$.  

As usual, $|D|$ denotes the support of the divisor $D$ of $\Sigma$. 
Then, for any $D$, there exists a
unique decomposition $D= D_{e}+ D_{s}$ such   that $|D_{e}|
\subset E$ and $\dim(|D_{s}|\cap E)<1$. 
Notice that $f\in m_{\mathcal{S},0}$ defines  an isolated
singularity at 0 if and only if 
$\mathrm{div}(f\circ p)_s$ in $\Sigma$ 
is reduced. The next theorem guarantees the existence of a function germ
$f$ with prescribed exceptional part $\mathrm{div}(f\circ p)_e$ which
is resolved by $p$, 
that is, such that $\mathrm{div}(f\circ p)$ is a normal crossing divisor. 
Notice that in such a case, the number $n_i$  of  components 
of $\mathrm{div}(f\circ p)_s$ (all of them smooth) intersecting $E_i$ is
determined by $D:=\mathrm{div}(f\circ p)_e$: 
indeed, $n_i=-D\cdot E_i$ for any $i$, as $\mathrm{div}(f\circ
p)\cdot E_i=0$.  

\begin{theorem}
\label{realisation} 
Let $p:(\Sigma,E)\to(\mathcal{S},0)$ be a good resolution of a normal
surface singularity  $(\mathcal{S},0)$ as above. Assume that the 
effective  divisor $D=\sum m_iE_i\neq 0$ satisfies 
$$ (D+E+K_{\Sigma})\cdot E_{i}+2 \leq 0 \ \ \ \ \mbox{for any
$i\in\{1,\ldots,r\}$}.$$
Then there exists a function $f\in m_{\mathcal{S},0}$, with an
isolated singularity at $0$, such that $\mathrm{div}(f\circ p)$ is a normal
crossing divisor on $\Sigma$   with $\mathrm{div}(f\circ p)_e=D$.
Moreover, for each $i$,
the number of intersection points $n_i=\mathrm{div}(f\circ p)_s\cdot E_i$
is strictly positive.
\end{theorem}

\begin{proof}
We refer to \cite{R 97} for an introduction to the methods used in this proof.

First notice that
\begin{equation*}
n_i=-D\cdot E_i\geq v_i+E_i^2+K_\Sigma\cdot E_i+2=v_i+2g_i.
\tag{$*$}\end{equation*}
The right member is strictly positive, except when $E=E_1$ is
just a rational curve. But in this case, $D=m_1E_1$ and $E_1^2<0$,
proving that $n_1>0$.
This answers the last statement.

For the existence of $f$, consider the exact sequence
$$0\to {\mathcal O}_\Sigma (-D-E)\to {\mathcal O}_\Sigma(-D)\to
{\mathcal O}_E(-D)\to 0.$$
Clearly, $H^1({\mathcal O}_\Sigma(-D-E))=
H^1({\mathcal O}_\Sigma(K_\Sigma+A))$, where $A:=-D-E-K_\Sigma$, with
$A\cdot E_i\geq 2$ for any $i$. Therefore, this last group is vanishing by the 
Laufer-Ramanujam theorem (see \cite{L 72} (3.2) or \cite{R 72}; it 
 is also called the `generalized' Grauert-Riemenschneider vanishing
theorem). This shows that $\pi:H^0({\mathcal O}_\Sigma(-D))\to
H^0({\mathcal O}_E(-D))$ is onto.  

Next, we show that ${\mathcal O}_E(-D)=\omega_E(A|_E)$ is generated by global 
sections. Firstly, if $q$ is a smooth point of $E$, then one has to show that
there exists a global section $s$ of $\omega_E(A|_E)$ with $s(q)$ nonzero.
For this it is enough to verify that $H^0(\omega_E(A|_E))\to
\omega_E(A|_E) \otimes {\mathcal O}_q$ is onto, or
$H^1(\omega_E(A|_E-q))=H^0({\mathcal O}_E(-A|_E+q))$ is zero.
But this last group is vanishing indeed,
since $-A+q$ is still negative on each component $E_i$. 
Similarly, if $q\in E_i\cap E_j$, with local coordinates $(x,y)$
such that $(E,q)=\{xy=0\}$, then it is convenient to consider the ideal
sheaf $I$ generated by $x+y$. Then $\omega_E(A|_E)\otimes I$ is still
locally free and its degree drops by one on both $E_i$ and $E_j$, hence
again $H^1(\omega_E(A|_E)\otimes I)=0$. 

Finally notice that the projective morphism 
$\varphi:E\to \PP^N$ induced by the globally generated 
${\mathcal O}_E(-D)$ is finite. Indeed, assume that $\varphi(E_i)$
is a point. Then the general hyperplane section of $\PP^N$
misses this point, which implies  $D\cdot E_i=0$ 
contradicting ($*$). 

Hence, if one takes $s\in H^0({\mathcal O}_\Sigma(-D))$ with
$\pi(s)$ a generic element of $H^0({\mathcal O}_E(-D))$,
then $s=f\circ p$  for some $f$ with the wanted properties.
\end{proof}

\begin{remark} \label{symmetry} A divisor $D$ which verifies the hypothesis
in this theorem always exists (since the intersection form is
negative definite). Moreover, one can also assume that
$D$ is fixed by the automorphism group of the weighted dual graph of $p$.
\end{remark}

\subsection*{B. Vertical links and horizontal open books in plumbed manifolds.}
Let $\Gamma$ be a finite connected weighted (plumbing) 
graph  with $r$ vertices 
$A_{1},...,A_{r}$. Each vertex $A_{i}$ is
weighted with two integers $(g_i,e_i)$, with $g_i\geq 0$.
Let $M(\Gamma)$ be the
oriented closed 3-manifold obtained from $\Gamma$ by plumbing (see
Mumford \cite{M 61} and Neumann \cite{N 81}). Briefly,  this 
construction runs as follows.
We associate with  each vertex $A_i$ an oriented circle
bundle $p_i :M_i \rightarrow S_i$ with Euler number $e_i$, 
where $S_i$ is an oriented compact
connected real surface of genus $g_i$,
then we glue these 3-manifolds according to the edges of the graph.

\begin{definition}
Let $M(\Gamma)$ be a plumbed 3--manifold. To any $r$-tuple
$\underline{n}=(n_1,\ldots,n_r)$ of non-negative integers 
one associates a link $N(\underline{n})$ in $M(\Gamma)$ 
as follows.  For each $i$,  consider $n_i$ generic fiber of the
 circle bundle $M_i\to S_i$, then $N(\underline{n})$ is their union, where 
 $i\in\{1,\ldots,r\}$. Any such link is called {\em vertical}.
\end{definition}

Any vertical link is naturally oriented  by the orientations of the fibers.

\begin{definition}
Let $M(\Gamma)$ be a plumbed manifold. A {\em horizontal open book}
in $M(\Gamma)$
is an open book $(N(\underline{n}),\theta)$ whose binding is a vertical link,
whose (open) pages are transversal to the fibers of the bundles $M_i\to S_i$,
and which is compatible with the orientations.
\end{definition}

\begin{example} 
The main example is provided  by normal surface
singularities and germs of analytic functions (see  \cite{N 81}).
Indeed, for any fixed good resolution $p$,
the {\em dual resolution graph} $\Gamma(p)$ serves as a plumbing graph, and 
the boundary $M(\mathcal{S})$ is diffeomorphic
to the plumbed manifold $M(\Gamma(p))$.
Moreover, if $f\in m_{\mathcal{S},0}$ defines an isolated singularity at $0$,
and $\mathrm{div}(f\circ p)$ is a normal crossing divisor,
then the Milnor open book $(N(f),\theta(f))$ defined
in $M(\mathcal{S})$ is isomorphic to a horizontal open book
$(N(\underline{n}),\theta)$
in $M(\Gamma(p))$, where each $n_i$ is exactly the number of
components of the strict transform $\mathrm{div}(f\circ p)_s$ which cut $E_i$.
\end{example}

\begin{proposition}\label{invariant}
Let $\Gamma$ be a connected weighted plumbing  graph with $r$ vertices whose
associated intersection form $I(\Gamma)$ is   non-degenerate.
Let $\underline{n}=(n_1,\ldots, n_r)$ be a $r$-tuple of \emph{strictly
  positive} integers.
Then any two horizontal open books in $M(\Gamma)$ with binding
$N(\underline{n})$ are
isomorphic.
\end{proposition}

\begin{proof}
In order to prove the proposition, we only have to collect  some existing
results from the literature (all of them stated in  \cite{P 01}).
The sequence of arguments is the following.

Any isomorphism class of open book can be characterized  completely by the
conjugacy class of the monodromy $h$ acting on the page $F$ of the open book. 
In the case of plumbed (or graph) manifolds, one can take for the  monodromy 
a quasi-periodical homeomorphism. In \cite{N 37} and \cite{N 44}
Nielsen associated with such a homeomorphism the (so called)
Nielsen graph. From this graph, in general, one cannot recover the conjugacy
class of $h$; but Chaves in \cite{C 96} completed this graph by some 
additional decorations -- in this way constructing the  `completed
Nielsen graph' -- and proved that this completed graph characterizes
completely the conjugacy class of the quasi-periodical homemomorphism $h$. 

On the other hand, Pichon in \cite{P 01} (see especially section 4)
describes the combinatorial relationship connecting the completed 
Nielsen graph and the plumbing graph $\Gamma'$ of the vertical link
$N(\underline{n})\subset M$. 
(Here $\Gamma'$ is obtained from the plumbing graph $\Gamma$ of $M$
by adding arrowheads corresponding to the link components.)
In general, the completed Nielsen graph 
contains more information, and cannot be recovered from $\Gamma'$
(see \ref{nonisot},b). Nevertheless, if each $n_i$ is strictly positive, then
$\Gamma'$ determines completely the completed Nielsen graph
(see Algorithm 4.8 of \cite{P 01}). 
This ends our proof as well. 

For the convenience of the reader we provide a few more details.
Recall that $\Gamma$ codifies the intersection matrix $I(\Gamma)=\{E_i\cdot
E_j\}_{i,j}$ represented in a fixed basis $\{E_i\}_i$. 
In the proof of \ref{realisation}, ($*$) provides $n_i$ as 
$-E_i\cdot (\sum_jm_jE_j)$. These  `multiplicities' $\{m_j\}_j$ 
constitute a part of the decorations of the Nielsen graphs. 
Since $I(\Gamma)$ is 
non-degenerate, they can be computed from the integers $\{n_i\}_i$. 

The plumbing construction provides a decomposition $M(\Gamma)=\cup_i
V_i$ of $M(\Gamma)$, where $V_i$ is an $S^1$-bundle over $S_i\setminus (v_i\ 
\mbox{discs})$. Let $F$ be the page of an open book with binding 
$N(\underline{n})$, and set $F_i:=F\cap V_i$ for each $i$. Let $r_i$ be the
number of connected components of $F_i$. It turns out that $r_i$ divides $m_i$.
The point is that, basically, the additional information contained
in the Nielsen graph (associated with the monodromy of the open book), 
compared with $\Gamma'$, is exactly the collection of
the integers $\{r_i\}_i$. It may help to think about this in the following
intuitive picture: one may construct a `covering graph' of $\Gamma'$ 
(which is equivalent with  the Nielsen graph) by
putting $r_i$ vertices above the vertex $A_i$ of $\Gamma'$ (and there is a
similar covering procedure for edges as well). If the integers $r_i$ are 
larger than one, it can happen that for the same collection of $\{r_i\}_i$'s
more covering types can appear. This global twisting data is the additional
information in the completed Nielsen graph, but which is superfluous
as soon as $r_i=1$ for all $i$. 

In our case, since $n_i>0$ for all $i$, analysing  tubular neighbourhoods
of the link components, we easily realize that each $F_i$ is connected,
i.e. $r_i=1$. Therefore, the completed Nielsen graph and $\Gamma'$ 
codify the same amount of information.
\end{proof}

\begin{remark}\label{nonisot}
 (a) In the context of \ref{invariant}, the two
 open books are not \emph{isotopic} in general. Indeed, one can take
 the image of any horizontal open book by a convenient 3-dimensional
 Dehn-twist, that is a self-diffeomorphism of $M(\Gamma)$ which is the
 identity outside a regular neighborhood of a``vertical torus" and
 which preserves the foliation by circles inside this regular
 neighborhood. 
\\
(b) If some components of $\underline{n}$ are allowed to vanish, then
 \ref{invariant} fails (see e.g. \cite{N 00} and \cite{P 01}). \\
(c) In \cite{CP 04}, the main result \ref{mthm} was established only
for Milnor fillable 3-manifolds which are rational homology spheres
because the authors 
 were not conscious of \ref{invariant}. Instead,
they  used the fact, easily deducible from results of Stallings and
Waldhausen, that in such a 3-manifold an open book is determined 
\textit{up to isotopy} by its binding alone.
\end{remark}

\begin{corollary}\label{ubique}
Let $M$ be a closed connected oriented 3-manifold which is
Milnor fillable. Then there exists an open book $(N,\theta)$
in $M$, which can be completely characterized by the topology of 
$M$,  such that, for any germ $(\mathcal{S},0)$ of normal complex
surface with
$M\simeq M(\mathcal{S})$, there exists a function
$f\in m_{\mathcal{S},0}$
having an isolated singularity at 0 whose Milnor open book
$(N(f),\theta(f))$ is
isomorphic to $(N,\theta)$.
\end{corollary}

\begin{proof} First notice that by the work of Neumann \cite{N 81} 
the homeomorphism type of $M$
determines the dual weighted graph
$\Gamma(p_0)$ of the minimal good resolution
$p_0:(\Sigma_0,E_0)\to(\mathcal{S},0)$
of \emph{any} Milnor filling $(\mathcal{S},0)$ of $M$.

Fix any of these fillings. 
Then choose an exceptional divisor $D$ in
$\Sigma_0$ as in  \ref{realisation}.
This choice only depends on $\Gamma(p_0)$, that is, on the topology of $M$.
Then apply \ref{realisation} and \ref{invariant}.
\end{proof}

\begin{remark}
 The notion of \textit{good resolution} used in \cite{N 81} does not
 coincide with the one we use here. The difference is that we ask the
 components $E_i$ to be smooth, while in \cite{N
 81} self-intersections are allowed. Nevertheless, minimal good resolutions
 exist and are unique for both 
 definitions, and the way to pass from one to the other can be
 completely described using the associated weighted dual resolution
 graphs. Hence, Neumann's result quoted before and proved in \cite{N 81}
 for his definition implies the analogous result for our definition. 
\end{remark}

\begin{proof}[{\bf C. Proof of Theorem \ref{mthm}.}]
Let $M$ be a closed connected oriented 3-manifold. By corollary \ref{ubique},
there exists an open book $(N,\theta)$ in $M$ which is isomorphic to a
Milnor open book
$(N(f),\theta(f))$ in any
Milnor filling $M(\mathcal{S})$ of $M$.
Now each of these Milnor open books carries the corresponding contact
structure
$(M(\mathcal{S}),\xi(\mathcal{S}))$ by theorem \ref{carry}. Since
they are all
isomorphic to $(N,\theta)$, theorem \ref{giroux}
shows that these contact structures are also all isomorphic.
\end{proof}

\section{Concluding remarks}
\label{persp}

In 3-manifolds which are circle bundles over oriented surfaces
of genus $>0$
with negative Euler number, 
there is only one isotopy class of contact structure transversal to the
fibers: see Giroux \cite[3.1 b), in conjunction with 3.1 a), 2.4 c)
and 2.5]{G 01}.  
Since any Milnor fillable contact structure on
such a manifold has this property, this shows that, in this case,
its \emph{isotopy} class is well defined (and not only its
\emph{isomorphism} class; 
we use here the fact, due to Waldhausen, that in this case the fibration
is well defined up to isotopy).
This, together with Remark \ref{nonisot} c), has to be compared with
Remark \ref{nonisot} a), where
we insisted on the (\emph{a priori}) non-isotopy  between horizontal open books
sharing the same binding.

Notice also that the homotopy type (as an \textit{unoriented}
2-plane field) of a Milnor fillable  contact 
structure is well defined (it is invariant up to isotopy by the group
$\mathit{Diff}^+(M)$ 
of self-diffeomorphisms  of $M$ which preserve the orientation). Firstly, two 
oriented plane fields which are positively transversal to the oriented
fibers of the plumbing structure are homotopic (as one can see by
taking an auxiliary Riemannian metric and rotating at each point one
plane into the other by the \textit{unique} shortest path). Secondly,
one can deduce from Neumann's work \cite{N 81} that the plumbing
decomposition (with \textit{unoriented fibres}) which corresponds to
the minimal good resolution is 
unique up to isotopy, which shows that it is invariant by
$\mathit{Diff}^+(M)$, up to isotopy. 

Moreover, using this last invariance and the results of \cite{CP 04},
we see that on 
a Milnor fillable rational homology sphere all Milnor fillable contact
structures are isotopic.

\smallskip

One can also ask the following general questions about the 
unique Milnor fillable contact
structure $\xi(M)$ (up to contactomorphism) on a Milnor fillable
3-manifold $M$. 

\smallskip

\textbf{Question 1} Characterize the subgroup of  $\mathit{Diff}^+(M)$ which
fixes $\xi(M)$ up to isotopy (here Remark \ref{symmetry} 
should be useful; as said before, we know that on Milnor fillable
rational homology spheres we get the whole group).

\textbf{Question 2} Find the Ozsv{\'a}th-Szab{\'o} contact invariant
associated with $\xi(M)$ (see \cite{OS 02} for its definition).


\textbf{Question 3} For each incompressible torus $T \subset M$, find
its associated torsion.

\smallskip

 The notion of \textit{torsion}
 of a contact structure originates in Giroux's article \cite{G 94}. 
It can be refined as a function defined on the set of isotopy
classes of incompressible tori in $M$, with values in $\N \cup
\{\infty\}$. See Colin, Giroux, Honda \cite{CGH 02} for its use in the
classification of tight contact structures on a given 3-manifold. 
If one believes in the manifestation in this context of a principle of
economy of algebraic geometry, one could conjecture that the torsion
of any incompressible torus in $M$ is equal to $0$.

\smallskip

\textbf{Question 4} Given a Milnor fillable 3-manifold, evaluate the
number of pairwise non-contactomorphic Stein fillable contact
structures on it. In particular, does there exist an example with an
infinite number of such structures?

\smallskip

Notice that on the lens spaces, which are Milnor fillable manifolds,
all tight contact structures are classified up to isotopy (see Giroux
\cite{G 00}): there is a finite number of them, they are all Stein
fillable and, 
with the exception of one or two, they are virtually overtwisted.

\end{document}